\documentstyle[amssymb,12pt]{article}

\begin{document}

\begin{center}
{\bf \Large Sigma Models, Minimal Surfaces and Some Ricci Flat Pseudo
Riemannian Geometries}\\

\vspace{2cm}

{\bf Metin G{\" u}rses} \\
{\small Department of Mathematics, Faculty of Science} \\
{\small Bilkent University, 06533 Ankara - Turkey}\\
{\it email:gurses@fen.bilkent.edu.tr}

\end{center}

\vspace{1cm}

\begin{center}
{\it
 Geometry, Integrability, and Quantization\\
July 7-15, 2000 Varna, Bulgaria\\
Eds. I.M. Mladenov and G.L. Naber}
\end{center}

\vspace{1cm}

\noindent
{\bf Abstract}:\,
We consider the sigma models where the base metric is proportional to the 
metric of the configuration space. We show that the corresponding
sigma model equation admits a Lax pair.
We also show that this type of sigma models in two dimensions
are intimately related to the minimal surfaces in a flat pseudo Riemannian
3-space. We define  two dimensional surfaces conformally 
related to the minimal surfaces in  flat three dimensional
geometries which enable us to give a construction of the metrics of some even  
dimensional Ricci flat (pseudo-) Riemannian geometries.

\section{Introduction}

Let $M$ be a 2-dimensional manifold with local coordinates
$x^{\mu}=(x,y)$
and $\Lambda^{\mu \nu}$ be the components of a tensor field in $M$. Let P be
an $2 \times 2$ matrix with a nonvanishimg constant determinant. 
We assume that $P$ is a hermitian ($P^{\dagger}=P$) matrix. Then the 
field equations of the sigma-model we consider is given as follows

\begin{equation}
{\frac{\partial}{\partial x^{\alpha}}}\,\Bigl(\Lambda^{\alpha \beta} P^{-1}
{\frac{\partial\, P}{\partial x^{\beta}}}\Bigr) = 0. \label{SM}
\end{equation}

\noindent
The integrability of the above equation has been studied in \cite{KAR}
where the matrix function $P$ and the tensor
$\Lambda^{\alpha\, \beta}$ were considered independent.
The sigma model equation given above is integrable provided
$\Lambda$ satisfies the conditions

\begin{equation}
\partial_{\alpha}\, ({1 \over \sigma}\, \Lambda^{\alpha \beta}\,
\partial_{\beta}\, \sigma)=0, ~~~ \partial_{\alpha}\, ({1 \over \sigma}\,
\Lambda^{ \beta \alpha}\, \partial_{\beta} \phi)=0
\end{equation}

\noindent
where $\sigma$ and $\phi$ are determinant and antisymmetric part of the
tensor $\Lambda^{\alpha \beta}$ respectively

We have classified in \cite{KAR}
possible forms of the tensor $\Lambda^{\alpha \, \beta}$ under these conditions
of integrability. The case where $\Lambda$ and $P$ are related has been 
considered in \cite{gur}.
As an example, let $P=g$ where $g$ is a $2 \times 2$ symmetric matrix .
Letting also $\Lambda^{\alpha \, \beta}=
g^{\alpha \, \beta}$ , the inverse components of the metric
 $g_{\alpha \, \beta}$ , then (\ref{SM}) becomes

\begin{equation}
{\frac{\partial}{\partial x^{\alpha}}}\,\Bigl(g^{\alpha \beta} g^{-1}
{\frac{\partial\, g}{\partial x^{\beta}}}\Bigr) = 0. \label{SM1}
\end{equation}

\noindent
The above sigma model equation is integrable and the Lax equation
is simply given by \cite{gur}

\begin{equation}
\epsilon^{\alpha\, \beta}\, \frac{\partial}{\partial x^{\beta}}\,\, \Psi=
{1 \over k^2+\sigma}\,(k\,g^{\alpha\, \beta}-\sigma\,\epsilon^{\alpha\, \beta})\,g^{-1}\,
\frac{\partial\,g}{\partial x^{\beta}}\,\, \Psi. \label{lp1}
\end{equation}

\noindent
Integrability conditions are satisfied because $det{g}=\sigma$ (a constant)
and $g$ is symmetric. Here $k$ is an arbitrary constant (the spectral
 parameter), $\epsilon^{\alpha \, \beta}$ is the Levi-Civita tensor with
$\epsilon^{12}=1$.

In the theory of surfaces in $R^{3}$ there is a class , the minimal surfaces,
which have special importance both in physics and mathematics \cite{bu1},
\cite{die}. Let $S=\{ (x,y,z) \varepsilon R^{3}; z=h(x,y)\}$ define a surface
$S \varepsilon R^3$ which is the graph of a differentiable function 
$\phi(x,y)$. This surface is called minimal if $\phi$ satisfies the condition

\begin{equation}
(1+\phi,_{x}^2)\, \phi,_{yy}-2 \phi,_{x}\, \phi,_{y}\, \phi,_{xy}+
(1+\phi,_{y}^2)\, \phi,_{xx}=0.
\label{met5}
\end{equation}

\noindent
The Gaussian curvature $K$ of the surface $S$ is given by

\begin{equation}
K={\phi,_{xx}\, \phi,_{yy}-\phi,_{xy}^2 \over (1+\phi,_{x}^2+\phi,_{y}^2)^2}.
\end{equation}

Here in this work we generalize the above treatment to more
general geometries. Instead of $R^3$ we take a pseudo Euclidean
manifold $M_{3}$ and two surfaces with any signature.

Let $(S, g)$ denote a two dimensional geometry where $S$ is a surface
in a three dimensional flat manifold $M_{3}$ and $g$  is a (pseudo-)
Riemannian metric on $S$ with a non vanishing determinant, $det(g)$.
Furthermore we assume that the metric components 
$g_{\alpha \beta}$ satisfies the following conditions

\begin{eqnarray}
\partial_{\mu}\, (g^{\mu \nu}\, g^{-1} \partial_{\nu}\,g)=0,\\
R+{1 \over 4}\, tr [g^{\mu \nu} \partial_{\mu}\, g^{-1} \partial_{\nu}\,g]=0,
\end{eqnarray}

\noindent
where $R$ is the Ricci scalar of $S$ We shall see in the following 
sections that
some surfaces which are conformally related to minimal surfaces  satisfy
the above conditions.

The importance of such surfaces arises when we are interested in
even dimensional Ricci flat geometries. By the utility the metric $g$ of
these surfaces we shall give a  construction (without solving any  further
differential equations) of the  metric of a $2N$ dimensional 
Ricci flat (pseudo-) Riemannian geometries. Ricci flat geometries are
important not only in differential geometry and general relativity
but also in gravitational instantons and in brane solutions of 
string theory \cite{brec}.

\section{Locally Conformal minimal surfaces}

Let $\phi$ be a differentiable function of $x^{1}=x$ and $x^{2}=y$ and $S_{0}$
be the surface in a three dimensional manifold $M_{3}$ 
 with a pseudo-Euclidean metric $g_{3}$ defined through
$ds^2=g_{0\, \mu \nu}\,dx^{\mu}\, dx^{\nu}+\epsilon \, (dx^{3})^{2}$, where
$\mu,\nu=1,2$ , $\epsilon= \pm 1$ and $g_{0}$ is a constant everywhere in 
$M_{3}$, invertible,
symmetric $2 \times 2$ matrix. In this work we assume Einstein summation
convention, i.e., the repeated indices are summed up. Let $S_{0}$ be  given as 
the graph of the function $\phi$, i.e., $S_{0}=\{ (x^{1},x^{2},x^{3})
\in M_{3} | x^{3}= \phi (x^{1},x^{2}) \}$. Then the metric on $S_{0}$ 
is given by

\begin{equation}
h_{\mu \nu}=g_{0\, \mu \nu}+\epsilon\, \phi_{, \mu}\, \phi_{, \nu}.
\label{2met}
\end{equation}

\noindent
Since\, $det\, h = (det \, g_{0})\, \rho$ where

\begin{equation}
\rho=1+\epsilon\, g_{0}^{\mu \nu}\, \phi_{, \mu}\, \phi_{,\nu}
\end{equation}

\noindent
then $h$ is everywhere (except at those points where $\rho=0$)
invertible. Its inverse is given by

\begin{equation}
h^{\mu \nu}=g_{0}\,^{\mu \nu}-{\epsilon \over \rho}\, \phi^{\mu}\, \phi^{\nu}
\end{equation}

\noindent
where $g_{0}^{\mu \nu}$ are the components of the inverse matrix $g_{0}^{-1}$
of $g_{0}$. Here the indices are lowered and raised by the metric $g_{0}$
and its inverse $g_{0}\, ^{-1}$ respectively. For instance , 
$\phi^{\mu}\,_{,\nu}=g_{0}^{\mu \alpha}\, \phi_{, \alpha \nu}$.
The Ricci tensor corresponding to the metric in  (\ref{2met}) is given by

\begin{equation}
r_{\mu \nu}={\epsilon \over \rho}\, (\nabla^{2}\, \phi)\, \phi_{,\mu \nu}
-{\epsilon \over \rho}\, \phi_{,\mu}\,^{\alpha}\, \phi_{,\nu \alpha}+
{1 \over 4 \rho^{2}}\, \rho_{,\mu}\, \rho_{, \nu}, \label{2ric}
\end{equation}

\noindent
where

\begin{equation}
\nabla^{2}\, \phi= h^{\mu \nu}\, \phi_{, \mu \nu}=
g_{0}\,^{\mu \nu}\, \phi_{, \mu \nu}-{ 1 \over  2 \rho}\, \phi^{\alpha}\,
\rho_{, \alpha}
\end{equation}

\noindent
The Ricci scalar or the Gaussian curvature $K$ and the mean curvature
$H$ are obtained as

\begin{eqnarray}
K&=&{\epsilon \over \rho^2}\, [(\phi^{\alpha}\,_{\alpha})^2-\phi^{\alpha \beta}\,
\phi_{,\alpha \beta}]\\
H&=&{1 \over \sqrt{\rho}}\, h^{\mu \nu}\, \phi_{,\mu \nu}, \label{min}
\end{eqnarray}

\noindent
The following equation is valid only for the case of two dimensional
geometries.

\begin{equation}
\phi_{,\alpha \mu}\, \phi_{,\beta \gamma}-\phi_{,\alpha \beta}\,
\phi_{,\mu \gamma}=
-\lambda_{0}\,(g_{0\, \alpha \mu}\,g_{0\,\beta \gamma}-
g_{0\,\alpha \beta}\,g_{0\,\gamma \mu}), \label{id1}
\end{equation}

\noindent
where

\begin{equation}
\lambda_{0}={1 \over 2}\,[\phi^{\alpha \beta}\, \phi_{\alpha \beta}-
(\phi^{\alpha}_{\alpha})^2].
\end{equation}

\vspace{0.3cm}

\noindent
Contracting this equation with $g^{\alpha \beta}$ leads to

\vspace{0.3cm}

$$\phi^{\alpha}_{\mu}\, \phi_{,\alpha \nu}-
\phi^{\alpha}_{\alpha}\, \phi_{,\mu \nu}=\lambda_{0}\,g_{0\,\mu \nu}$$

\noindent
We also have

\vspace{0.3cm}

$$r_{\alpha \beta}={K \over 2}\, h_{\alpha \beta},~~~\lambda_{0}=
-{\epsilon \over 2}\,\rho^2\,K$$

\vspace{0.3cm}

\noindent
For the minimal surfaces we have $H=0$ and the following important
properties of the metric $h_{\alpha \beta}$ on $S_{0} \cite{gur1}$

\vspace{0.3cm}

\begin{eqnarray}
\partial_{\alpha}\,[\sqrt{\rho}\, h^{\alpha \beta}\,
\partial_{\beta}\, \phi]&=&0, \label{sig1}\\
\partial_{\alpha}\,( \sqrt{\rho}\, h^{\alpha \beta})&=&0. \label{sig2}
\end{eqnarray}

\vspace{0.3cm}

\noindent
We now define surfaces which are locally conformal 
to minimal surfaces. Let $S$ be such a surface, i.e., locally conformal
to $S_{0}$. Then the  metric on $S$ is given by

\begin{equation}
g_{\alpha \beta}={ 1 \over \sqrt{\rho}}\, h_{\alpha \beta}.
\end{equation}

\noindent
It is clear that $det\, g= det \,g_{0} \ne 0$. In the sequel we shall assume
that the surface $S_{0}$ is minimal and hence the metric defined on it
satisfies all the equivalent conditions in (\ref{sig1}) and (\ref{sig2}). 
The corresponding Ricci tensor of $g$ is given as

\begin{equation}
R_{\alpha \beta}=r_{\alpha \beta}-(\nabla_{g}^{2}\, {\psi_{0}})\,
g_{\alpha \beta},
\end{equation}

\noindent
where $\psi_{0}=-{1 \over 4}\, log (\rho)$ and $\nabla_{g}^{2}$ is 
the Laplace-Beltrami operator with respect to
the metric $g$. We then have

\vspace{0.3cm}

\noindent
{\it Proposition 1.\, The following equation is an identity related
to the conformal surface $S$.

\begin{equation}
R=-{1 \over 4} g^{\alpha \beta}\, tr [\partial_{\alpha}\, g^{-1}\, 
\partial_{\beta}\, g].
\end{equation}

\noindent
Here $g$ is the $2 \times 2$ matrix of $g_{\alpha \beta}$ and $g^{-1}$
is its inverse. The operation tr is the standard trace operation for
matrices.}

\vspace{0.3cm}

\noindent
In the following parts of the work we need some harmonic functions
with respect to the metric $g$. For this purpose  we introduce some 
vectors on $S$.
Let $v_{\alpha}=(1,0),~v^{\prime}_{\alpha}=(0,1)$ and
$u^{\alpha}=(1,0),~u^{\prime\, \alpha}=(0,1)$. We now define some functions
over $S$.

\begin{eqnarray}
\xi_{1}&=&g^{\alpha \beta}\, v_{\alpha}\, v_{\beta},~~~\xi_{2}=
g^{\alpha \beta}\, v^{\prime}_{\alpha}\, v^{\prime}_{\beta},\\
w_{1}&=&\sqrt{\rho}\, g_{\alpha \beta}\,u^{\alpha}\, u^{\beta},~~~
w_{2}=\sqrt{\rho}\, g_{\alpha \beta}\, u^{\prime \, \alpha}
\,u^{\prime\, \beta}.
\end{eqnarray}

\noindent
It is now easy to prove

\vspace{0.3cm}

\noindent
{\it Proposition 2.

\begin{eqnarray}
\nabla_{g}\,^2 \zeta -a_{0}\, R&=&-a_{0}\, \sqrt{\rho}\,K ,
\label{oz01}\\
\nabla_{g}\,^2 \psi_{1}-(a_{1}+a_{2})\, R&=&0,  \label{oz02}\\
\nabla_{g}\,^2 \psi_{2}-2(b_{1}+b_{2})\, R&=&-(b_{1}+b_{2})\,\sqrt{\rho}\,K,
\end{eqnarray}

\noindent
where

\begin{eqnarray}
\zeta&=&{a_{0} \over 2}\, log (\rho),\\
\psi_{1}&=&a_{1}\,log(\xi_{1})+a_{2}\, log(\xi_{2}),\\
\psi_{2}&=&b_{1}\, log(w_{1})+b_{2}\, log(w_{2}).
\end{eqnarray}

\noindent
Here $a_{0},a_{1},a_{2},b_{1},$ and $b_{2}$ are arbitrary constants.}

\vspace{0.3cm}

\noindent
The function $\mu$ defined by $\mu=(b_{1}+b_{2})\, \zeta-a_{0}\, \psi_{2}$
satisfies similar equation as $\psi_{1}$

\begin{equation}
\nabla_{g}\,^{2}\, \mu=-a_{0}\,(b_{1}+b_{2})\,R, \label{mu}.
\end{equation}

\noindent
Hence we have two different solutions of the equation

\begin{equation}
\nabla_{g}^{2} \sigma=-{c \over 4}\, g^{\alpha \beta}\, tr [
(\partial_{\alpha}\,g^{-1})\, \partial_{\beta}\, g], \label{cc1}
\end{equation}

\noindent
for some function $\sigma$. If $\sigma=\psi_{1}$ then $c=a_{1}+a_{2}$,
if $\sigma=\mu$ then $c=-a_{0}(b_{1}+b_{2})$. It is straightforward 
 to show that 

\begin{equation}
\xi_{1}={w_{2} \over det\, g_{0}\, \sqrt{\rho}},~~ 
\xi_{2}={w_{1} \over det\, g_{0}\,\sqrt{\rho}}
\end{equation}

\noindent
Hence $\psi_{1}$ will not be considered as an independent function.
It is  interesting and important to note that under the minimality
condition the matrix $g$ satisfies the following condition as well.

\vspace{0.3cm}

\noindent
{\it Proposition 3.\, Minimality of $S_{0}$, $H=0$, also implies
a sigma model \cite{gur1}, \cite{gur2} like equation for $g$, i.e.,

\begin{equation}
\partial_{\alpha}\,[g^{\alpha \beta}\, g^{-1}\,
\partial_{\beta}\,g]=0. \label{sig}
\end{equation}
}

\vspace{0.3cm}

\noindent
{\it Proof}:\,  The metric $g_{\alpha \beta}$ and
its inverse $g^{\alpha \beta}$ are written in a nice  form

\begin{eqnarray}
g_{\alpha \beta}={1 \over \sqrt{\rho}}\,(g_{0\,\alpha \beta}+\epsilon
\phi_{,\alpha}\, \phi_{,\beta}),\\
g^{\alpha \beta}=\sqrt{\rho}\, (g_{0}^{\alpha \beta}-
{\epsilon  \over \rho}\, \phi^{\alpha}\, \phi^{\beta})
\end{eqnarray}

\noindent
where $g_{0\,\alpha \beta}$ are the components of the  matrix $g_{0}$.
The minimality condition $H=0$ reduces to
$g^{\alpha \beta}\, \phi_{,\alpha \beta}=0$ or

\begin{equation}
\phi^{\alpha}\,_{,\alpha}={\phi^{\alpha}\, \rho_{\alpha} \over 2 \rho}.
\label{min2}
\end{equation}

\noindent
This condition also implies  

\begin{equation}
\partial_{\mu}\, g^{\mu \nu}=0. \label{min3}
\end{equation}

\noindent
Hence the sigma model equation (\ref{sig}) to be proved takes the form

\begin{equation}
h^{\mu \nu}\, \partial_{\nu}\,[g^{\alpha \gamma}\, \partial_{\mu}\,
g_{\gamma \beta}]=0, \label{sig0}
\end{equation} 

\noindent
where $h_{\alpha \beta}=\sqrt{\rho}\, g_{\alpha \beta}$.
It is straightforward to show that

\begin{eqnarray}
(g^{-1}\, \partial_{\mu}\,g)^{\alpha}_{\beta}&=&g^{\alpha \gamma}\,
\partial_{\mu}\, g_{\beta \gamma} \nonumber \\
&=&-{1 \over 2} {\rho_{,\mu}
\over \rho}\,\delta^{\alpha}_{\beta}-
{\epsilon \over 2}\, {\rho_{,\mu} \over \rho}\,
\phi^{\alpha}\, \phi_{,\beta}+\epsilon\phi_{,\beta}\, \phi^{\alpha}\,_{, \mu}+
{ \epsilon \over \rho}\, \phi^{\alpha}\, \phi_{, \mu \beta}, \label{den1}
\end{eqnarray}

\noindent
Using the identity (\ref{id1}) and the minimality condition 
(\ref{min2}) we obtain the following

\begin{eqnarray}
\rho_{,\mu}\, \phi_{, \beta \gamma}-\rho_{, \beta}\, \phi_{, \mu \gamma}=
2\epsilon \,\lambda_{0}\,(\phi_{, \mu}\, g_{0\, \beta \gamma}- 
\phi_{, \beta}\, g_{0\, \gamma \mu}),\\
\rho^{\mu}\, \phi_{, \mu \beta}=\phi^{\alpha}\,_{, \alpha}\, \rho_{, \beta}-
2 \epsilon\,\lambda_{0}\,\phi_{, \beta},\\
h^{\alpha \beta }\,\partial_{\alpha}\, ({ 1 \over \rho}\, 
\partial_{\beta}\, \rho)+{2\, \lambda_{0}  \over \rho^2}\,(1+\rho)=0.
\end{eqnarray}

\noindent
Utilizing these identities we get

\begin{eqnarray}
h^{\alpha \beta}\, \phi^{\mu}\,_{, \alpha \beta}=-
{2 \epsilon\,\lambda_{0} \over \rho}\,\phi^{\mu},\\
h^{\alpha \beta}\, \phi_{, \mu \alpha} \, \phi^{\nu}\,_{, \beta}=
-{\lambda_{0} \over \rho}\, \delta^{\nu}\, _{\mu}-
{\epsilon\,\lambda_{0} \over \rho}\,\phi^{\nu}\, \phi_{, \mu},\\
h^{\alpha \beta}\, \partial_{\alpha}\,(\phi^{\nu}_{\, , \beta}\,
\phi_{, \mu})=
-{\lambda_{0} \over \rho}\, \delta^{\nu}_{\mu}-
{3 \epsilon\,\lambda_{0} \over \rho}\,\phi^{\nu}\, \phi_{, \mu},\\
\rho_{, \alpha}\,h^{\alpha \beta}\, \partial_{\beta}\,(\phi^{\nu}\, 
\phi_{, \mu})
=-4 \epsilon\, \lambda_{0}\,\phi^{\nu}\, \phi_{, \mu}.
\end{eqnarray}

\noindent
Now applying $\partial_{\nu}$ to  (\ref{den1}) then multiplying by 
$h^{\mu \nu}$ and using the above identities (by virtue of the
minimality condition (\ref{min2}) ) it is easy to show (\ref{sig0}).

\noindent
Hence for every minimal surface $S_{0}$ and its metric $h$ we have a
conformally related surface $S$ with metric $g={h \over \sqrt{\rho}}$ 
($\det{h}=\rho\, \det{g_{0}}$) satisfying the conditions

\begin{eqnarray}
R+{1 \over 4}\, g^{\alpha \beta}\,tr [ \partial_{\alpha} \,g^{-1}\, 
\partial_{\beta}\,g]=0,\\
\partial_{\alpha}\,[g^{\alpha \beta}\,g^{-1}\, \partial_{\beta}\,g]=0.
\end{eqnarray}

\noindent
Here $g$ has determinant equals to $det{g_{0}}$ which is a nonzero constant.
This does not violate the covariance of our formulation because we could
formulate everything in terms of the metric $h$ of the minimal surfaces 
$S_{0}$ but the  above identities become lengthy and complicated. We loose 
no generality by using surfaces $S$ and the metric $g$ on them.

\section{Ricci flat pseudo Riemannian geometries}

We start first with four dimensions.
Let the metric of a four dimensional manifold $M_{4}$ be given by

\begin{equation}
ds^2=e^{2\psi}\,g_{\alpha \beta}\,dx^{\alpha}\,dx^{\beta}+
\epsilon_{1}\,g_{\alpha \beta}\,dy^{\alpha}\,dy^{\beta}, \label{4met}
\end{equation}

\noindent
where $\psi$ is a function of $x^{\alpha}$ and $\epsilon_{1}=\pm 1$.
Local coordinate of $M_{4}$ are denoted as $x^{a}=(x^{\alpha},
y^{\alpha}), ~a=1-4$

\vspace{0.3cm}

\noindent
{\it Proposition 4 .\, The Ricci flat equations $R_{ab}=0$ for the metric
(\ref{4met}) are given in two sets. One set satisfied identically
due to the Proposition 3 above and the second one is given by

\begin{equation}
\nabla_{g}\,^{2} \psi=0. \label{psi2}
\end{equation}
}

\vspace{0.3cm}

\noindent
There are two independent functions satisfying the above
Laplace equation , $\phi$ and $\mu$.
Using (\ref{mu}) we find that $\psi=e_{0}\, \phi+e_{1}\, \mu$
where $e_{0}$ and $e_{2}$ are arbitrary constants and $b_{2}=-b_{1}$.
Combining all these constants we find that

\begin{equation}
e^{2\psi}=e^{2 e_{0} \, \phi}\,w_{1}^{-2m_{1}}\,w_{2}^{-2m_{2}},
\label{psi}
\end{equation}

\noindent
where $m_{1}$ and $m_{2}$ are constants satisfying $m_{1}+m_{2}=0$.
Then the line element (\ref{4met}) becomes

\begin{equation}
ds^{2}={e^{2\,e_{0}\, \phi} \over w_{1}^{2m_{1}}\,w_{2}^{2m_{2}}} 
\,{h_{\alpha \beta}\,dx^{\alpha}\,dx^{\beta} \over \sqrt{\rho}}+
{h_{\alpha \beta}\,dy^{\alpha}\,dy^{\beta} \over \sqrt{\rho}}, \label{ins}
\end{equation}

\noindent
where $\phi$ satisfies the minimality condition ($H=0$)  (\ref{min})  which 
is explicitly given by

\begin{equation}
[k_{2}+\epsilon\,(\phi_{,y})^{2}]\,\phi_{,xx}-2[k_{0}+
\epsilon\, \phi_{,x}\, \phi_{,y}]\, \phi_{,xy}+[k_{1}+\epsilon\,
(\phi_{,x})^{2}]\, \phi_{,yy}=0, \label{min1}
\end{equation}

\noindent
where we take $(g_{0})_{ 11}=k_{1},~(g_{0})_{01}=k_{0},~(g_{0})_{22}=k_{2}$
and   assume that $det\, (g_{0})=k_{1}\,k_{2}-k_{0}^{2} \ne 0$. Hence
the functions $w_{1}$ and $w_{2}$ are given explicitly as

\begin{equation}
w_{1}=k_{1}+\epsilon\, (\phi_{,x})^{2},~~w_{2}=k_{2}+
\epsilon\, (\phi_{,y})^{2}. \label{w12}
\end{equation}

\noindent
The metric in (\ref{ins}) with  $e_{0}=0, m_{1}=m_{2}=0$ reduces to 
an instanton metric \cite{yn}. 

We shall now generalize Proposition 4 for an arbitrary even dimensional
pseudo-Riemannian geometry.
Let $M_{2+2n}$ be a $2+2n$ dimensional manifold with a metric

\begin{equation}
ds^2=e^{2\Phi}\,g_{\alpha \beta}\,dx^{\alpha}\,dx^{\beta}
+G_{AB}\,dy^{A}\,dy^{B}, \label{nmet}
\end{equation}

\noindent
where the local coordinates of $M_{2+2n}$ are given by $x^{\alpha+A}=
(x^{\alpha},y^{A}),~A=1,2, \cdots ,2n$, $\Phi$ and $G_{AB}$ are functions
of $x^{\alpha}$ alone. The Einstein  equations are given in
the following proposition

\vspace{0.3cm}

\noindent
{\it Proposition 5.\, The Ricci flat equations for the metric in
(\ref{nmet}) are given by

\begin{eqnarray}
\partial_{\alpha}\,[g^{\alpha \beta}\,G^{-1}\, \partial_{\beta}\,G]=0,\\
\nabla_{g}\,^{2}\, \Phi={1 \over 8}\, g^{\alpha \beta}\,
tr [(\partial _{\alpha} G^{-1}) \, \partial_{\beta}\,G]+{R \over 2},
\end{eqnarray}

\noindent
where $G$ is $2n \times 2n$ matrix of $G_{AB}$ and $G^{-1}$ is
its inverse.}

\vspace{0.3cm}

\noindent
Let us choose $G$ as a block diagonal matrix and
each block is the $2 \times 2$ matrix $g$. This means that
the metric in (\ref{nmet}) reduces to a special form

\begin{eqnarray}
ds^2=e^{2\Phi}\, g_{\alpha \beta}\,dx^{\alpha}\,dx^{\beta}+
\epsilon_{1}\,g_{\alpha \beta}\,dy_{1}^{\alpha}\,dy_{1}^{\beta}+
\cdots +\epsilon_{n}\,g_{\alpha \beta}\,dy_{n}^{\alpha}\,
dy_{n}^{\beta}, \label{smet}
\end{eqnarray}

\noindent
where the local coordinates of $M_{2+2n}$ are given by $x^{\alpha+A}=
(x^{\alpha},y_{1}^{\alpha}, \cdots , y_{n}^{\alpha})$, $\epsilon_{i}
=\pm 1,~ i=1,2, \cdots , n$. Then we have the following  theorem

\vspace{0.3cm}

\noindent
{\it Theorem.\, For every two dimensional minimal surface $S_{0}$
immersed in a three
dimensional manifold $M_{3}$ there corresponds a $2N=2+2n$-dimensional
Ricci flat (pseudo-) Riemannian geometry with the metric given in (\ref{smet})
with

\begin{equation}
e^{2\, \Phi}=e^{2\, \psi}\, \,{w_{1}^{-2n_{1}}\, w_{2}^{-2n_{2}}}
\,\, \rho^{n_{1}+n_{2} } , \label{phi}
\end{equation}

\noindent
where $\psi$ is given in (\ref{psi}), $w_{1}$ and $w_{2}$ are given
in (\ref{w12}), $n_{1}$ and $n_{2}$ satisfy

\begin{equation}
n_{1}+n_{2}={n-1 \over 2}. \label{n12}
\end{equation}
}

\vspace{0.3cm}

\noindent
{\it Proof:}\, Using proposition 5 for the metric
(\ref{smet}) the Ricci flat equations reduce to the following
equation

\begin{equation}
\nabla_{g}\,^{2}\, \Phi={n-1 \over 8}\, g^{\alpha \beta}\, tr
[ (\partial_{\alpha} g^{-1})\, \partial_{\beta} g]
\end{equation}

\noindent
By using (\ref{cc1}) and letting $a_{0}\,b_{1}=n_{1},~
a_{0}\,b_{2}=n_{2}$ and $\Phi=\mu+\psi$ we find (\ref{phi}) with the condition
(\ref{n12}). Here $\psi$ is a harmonic function ({\ref{psi2}) 
with respect to the metric
$g_{\alpha \beta}$. A solution of this function is given in the previous 
section in  (\ref{psi}).
Metric functions $\psi$, $w_{1}$, $w_{2}$ and $g_{\alpha \beta}$
are expressed explicitly in terms the function $\phi$ and its
derivatives $\phi_{,x}$ and $\phi_{,y}$. This means that for each solution
$\phi$ of (\ref{min1}) there exists a $2N$-dimensional metric (\ref{smet}).

The dimension of the manifold is $4\,(1+n_{1}+n_{2})$.
Here $n=1$ or $n_{1}+n_{2}=0$ corresponds to the four dimensional case. 
The signature of the geometry depends on the 
signature of $S$. If $S$ has zero signature then $M_{2N}$ has also
zero signature, but if the signature of $S$ is $2$ then the signature 
of $M_{2N}$ is $2\,(1+\epsilon_{1}+ \cdots +\epsilon_{n})$.

\vspace{1cm}

This work is partially supported by the Scientific and Technical
Research Council of Turkey (TUBITAK) and Turkish Academy of Sciences (TUBA).

\end{document}